\documentclass[oneside]{amsart}
\usepackage{amsmath,ifthen,amscd,amssymb,
graphicx,subfigure}
\newcommand{\BBMS}{MR98j:20034}
\newcommand{\BestBrady}{MR98i:20039}
\newcommand{\BB}{BB}
\newcommand{\BBIsop}{MR2001j:20046}
\newcommand{\Bowditch}{MR96b:20046}
\newcommand{\Bridson}{MR2000g:20071}
\newcommand{\BridsonHaefliger}{MR1744486}
\newcommand{\Bridcat}{MR1851271}
\newcommand{\Cannon}{MR88a:20049}

\newcommand{\DicksLeary}{MR99c:20050}
\newcommand{\Epatterns}{MR1974064}
\newcommand{\Ethesis}{Ethesis}
\newcommand{\Enonhopf}{MR2022477}
\newcommand{\Efftpft}{MR1950887}
\newcommand{\Epstein}{MR93i:20036}
\newcommand{\Gerstencat}{MR94j:20043}
\newcommand{\GerstenGGT}{MR96k:20073}

\newcommand{\NSgeomfinite}{MR96c:20066}
\newcommand{\NibloR}{MR99a:20037}
\newcommand{\Olshan}{MR93d:20067}

\newcommand{\Stallings}{MR28:2139}
\newcommand{\Thiel}{MR95e:20052}
\newcommand{\Wise}{MR96m:20058}

\graphicspath{{.}{./eps/}}

\newtheorem{thm}{Theorem}[section]
\newtheorem{lem}[thm]{Lemma}

\newtheorem{cor}[thm]{Corollary} 

\theoremstyle{definition}
\newtheorem{defn}[thm]{Definition}

\newtheorem{eg}{Example}

\newcommand{\blackboard}[1]{\ensuremath{\mathbb{#1}}}

\newcommand{\smallcaps}[1]{\textrm{\textsc{#1}}}

\newcommand{\N}{\blackboard{N}}
\newcommand{\Z}{\blackboard{Z}}

\newcommand{\R}{\blackboard{R}}

\newcommand{\cat}{\smallcaps{CAT}}

\newcommand{\ra}{\rightarrow}
\newcommand{\iif}{isoperimetric function}
\newcommand{\ac}{almost convex}
\newcommand{\fftp}{falsification by fellow traveler property}
\newcommand{\lsp}{loop shortening property}
\newcommand{\cg}{Cayley graph}
\newcommand{\hnn}{HNN extension}
\newcommand{\fg}{finitely generated}
\newcommand{\gset}{generating set}

\newcommand{\gam}{\gamma}
\newcommand{\ga}{$\gamma$}

\begin{document}

\title[The loop shortening property]
  {The loop shortening property and\\
  almost  convexity    }

\author[M.~Elder]{Murray J. Elder}
\address{(As at Dec 2006) Dept of Mathematical Sciences\\
	Stevens Institute of Technology\\
	Hoboken NJ 07030 USA}
\email{melder@stevens.edu}

\date{Published in  Geometriae Dedicata Vol 102 No. 1 (2003) 1--18}

\begin{abstract}
We introduce the \lsp\ and the basepoint \lsp\ for  \fg\ groups, and examine their relation
to quadratic \iif s and \ac ity.
\end{abstract}

\maketitle


\section{Introduction}
In this article we introduce two new properties of groups: the {\em \lsp}
and the {\em basepoint \lsp}.
The properties are  natural generalizations of the \fftp\ introduced by
Neumann and Shapiro \cite{\NSgeomfinite},
 which in turn is closely related to the property
of being \ac, introduced by Cannon in \cite{\Cannon}.
The first part of the article is devoted to proving three facts.
In Theorem  \ref{thm:asynch=synch} we see that asynchronous and
synchronous versions of
both properties are equivalent.
Theorem \ref{loopquad} states that having the \lsp\ implies 
finite presentability and 
a quadratic Dehn function, and Theorem
\ref{bploopac} shows that if a group presentation has the basepoint
\lsp\ then it is \ac.
In the second part of the article we examine four group presentations,
which  exhibit a diverse spectrum of properties. These examples
answer several natural questions about the loop shortening properties and the
interdependence between them and \ac ity, the \fftp\ and quadratic \iif s.

The author is indebted to
Noel Brady,  Jon McCammond and Walter Neumann for their
ideas and suggestions with this paper. 
In addition the author wishes
to thank an anonymous reviewer for 
her/his careful reading and suggestions.

\section{Preliminaries}


Throughout this article let $(G,X)$ denote the pair of a group and a finite
\gset. The set $X^*$ denotes the set of all words in the letters of $X$,
including the empty word. The \cg\ for the pair is denoted
$\Gamma(G,X)$.

\begin{defn}[Path,loop]
A word $w\in X^*$ corresponds to a {\em path} based at some vertex of the
\cg.
A {\em loop} is a path which starts and ends at the same vertex.
A path [loop] can be parameterized  by arc length, and we denote the
point at distance $t$ along the path [loop] from the start point by
$w(t)$.
For $t>|w|, w(t)$ is defined to be the endpoint of $w$.
\end{defn}

\begin{defn}[(Asynchronous) fellow traveling]
Two paths [loops] $w,u$ are said to {\em $k$-fellow travel} if 
$d(w(t),u(t))\leq k$ for all $t\geq 0$, where $w$ and $u$ can have
distinct start and end points.
Two paths [loops] $w,u$ are said to {\em asynchronously $k$-fellow
travel} if there is a proper monotone 
increasing  function $\phi:\R_{\geq 0} \ra \R_{\geq 0}$
such that 
$d(w(t),u(\phi(t))\leq k$ for all $t\geq 0$, where $w$ and $u$ can have
distinct start and end points.
\end{defn}

\begin{defn}[Falsification by fellow traveler property]
$(G,X)$ enjoys the [{\em asynchronous}] {\em \fftp}  if there
is a constant $k>0$ such that for each non-geodesic path $w$  in $(G,X)$,
 there is a path $u$ in $(G,X)$ with the same endpoints so that $|u|<|w|$ and $w,u$
[asynchronously] $k$-fellow travel.
\end{defn}

\noindent
Neumann and Shapiro introduced this property in \cite{\NSgeomfinite},
where they prove that if a pair $(G,X)$ enjoys the property, then the full
language of geodesics in the generators $X$ is regular.
They also show that the property is dependent on choice of \gset.
A reasonably simple proof in \cite{\Efftpft} shows that the
asynchronous and synchronous versions of this property are in fact equivalent.

\begin{defn}[Loop shortening property]
$(G,X)$ enjoys the {\em [asynchronous] loop shortening property} if there
is a constant $k>0$ such that for each loop $w$  in $(G,X)$,
 there is a loop $u$  in $(G,X)$ so that $|u|<|w|$ and $w,u$
[asynchronously] $k$-fellow travel.
\end{defn}

\noindent
Note that $u$ and $w$ can be disjoint.
A (seemingly) stronger version of the property is the following.

\begin{defn}[Basepoint loop shortening property]
$(G,X)$ enjoys the {\em [asynchronous] basepoint loop shortening property} if there
is a constant $k>0$ such that for each loop $w$  in $(G,X)$ based at
$w(0)$, there is a loop $u$  in $(G,X)$ based at $w(0)$ so that $|u|<|w|$ and $w,u$
[asynchronously] k-fellow travel.
\end{defn}

\begin{defn}[Almost convex]
$(G,X)$ is {\em almost convex} if there is a constant $C$ such
that every pair of points lying distance at most 2 apart and within
distance $N$ of the identity in $\Gamma(G,X)$ are connected by a path
of length at most $C$ which lies within distance $N$ of the
identity.
\end{defn}

\noindent
See \cite{\Cannon} for properties of \ac\ groups.
This property also depends on the choice of \gset\ \cite{\Thiel}.

Let $R\subseteq X^*$ denote some set of relators such that $(G,X)$
admits the presentation $\langle X,R\rangle$. Let $F(X)$ be the free
group generated by $X$, and let $N(R)$ be the normal closure of $R$ in
$F(X)$.  A word in $X^*$ represents the identity in $G$ if and only if
it is freely equal to an expression of the form
$$\Pi_{i=1}^k g_ir_ig_i^{-1}$$ where the $g_i\in F(X)$ and $r_i\in
R\cup R^{-1}$.  Define the  {\em area} $A(w)$ of a word $w\in X^*$
which represents the identity to be the minimum $k$ in any such
expression for $w$. 

\begin{defn}[Dehn function, Isoperimetric function]
 A {\em Dehn function} for $\langle X|R\rangle$ is defined to be
$\delta(n)=\max\{A(w):|w|\leq n\}$.  An {\em isoperimetric function}
for $\langle X|R\rangle$ is any function which satisfies $f(n)\geq
\delta(n)$.
\end{defn}

\noindent
Two functions $f,g$ are said to be equivalent if there are constants
$A,A',B,B',C,C',$

\noindent
$D, D',E,E'$ so that $f(n)\leq Ag(Bn+C)+Dn+E$, and
$g(n)\leq A'f(B'n+C')+D'n+E'$.  With respect to this definition, a
Dehn function of a group is \gset\ independent.  If $G$ has a
sub-quadratic \iif\ then its Dehn function is linear
\cite{\Bowditch,\Olshan}.  The class of groups which have a quadratic Dehn
function is diverse and not particularly well understood.  Examples
include \cat(0) groups \cite{\BridsonHaefliger}, automatic groups and
$(2n+1)$-dimensional integral Heisenberg groups for $n\geq 2$
\cite{\Epstein}.

\section{Asynchronous versus synchronous}

In this section we prove that the asynchronous and synchronous
 versions of the two properties are equivalent.

\begin{lem}[Discrete to continuous]\label{lem:discrete}
Let $w,u$ be paths in $(G,X)$, parameterized by $t\in \mathbb R_{\geq
0}$. If for some constant $k$ $d(w(t),u(t))\leq k$ for all $t\in
\mathbb N$ then $d(w(t),u(t))\leq k+1$ for all $t\in \mathbb R_{\geq
0}$.
\end{lem}

\begin{proof}
If $t-\lfloor t \rfloor \leq \frac{1}{2}$ then there is a path of
length at most $k+1$ from $w(t)$ to $w(\lfloor t \rfloor)$ to
$u(\lfloor t \rfloor )$ to $u(t)$. If $t-\lfloor t \rfloor >
\frac{1}{2}$ then there is a path of length at most $k+1$ from $w(t)$
to $w(\lceil t \rceil)$ to $u(\lceil t \rceil)$ to $u(t)$.
\end{proof}

\noindent
It follows that in order to prove that two paths synchronously
$k$-fellow travel it is sufficient to show that integer points are
within $(k-1)$ of each other.

\begin{thm}
\label{thm:asynch=synch}
 $(G,X)$ has the asynchronous [basepoint] \lsp\ if and only if $(G,X)$
  has the synchronous [basepoint] \lsp.
\end{thm}

\begin{proof}
If $(G,X)$ has the synchronous [basepoint] \lsp\ then it clearly has
the asynchronous [basepoint] \lsp.  Let $w$ be a loop of length $n$ in
$(G,X)$. If $(G,X)$ has the asynchronous \lsp\ with constant $k>0$
then there is a shorter loop $u$ and a proper monotone increasing
function $\phi:\R_{\geq 0} \rightarrow \R_{\geq 0}$ such that
$d(w(t),u(\phi(t)))\leq k$ for all $t\in \R_{\geq 0}$.  Without loss
of generality we may assume that $\phi(0)=0$.  Fix some constant
$\epsilon$ so that $0<\epsilon<1$.

Let $E=\{t\in \R_{\geq 0}:\phi(t)=|u| \}$. Note that if $t\in E$ then
$s\in E$ for all $s\geq t$.  If $n\notin E$ then define
$\phi':\R_{\geq 0} \rightarrow \R_{\geq 0}$ by
\[
\phi'(t)=  \{	
\begin{array}{lcr}
 \phi(t) & & 0\leq t \leq n-\epsilon \\
 \phi(n)+(|u|-\phi(n))\frac{t+\epsilon-n}{\epsilon} & &  n-\epsilon \leq t \leq n\\
 |u| & & t\geq n
\end{array}
\]

\noindent
Note that all points on $u$ from $u(\phi(n))$ to $u(|u|)$ are at most
$k$ from the vertex $w(n)$, as seen in Figure \ref{figasynch1}.
\begin{figure}[ht!]
  \begin{center}
     \includegraphics[width=5cm]{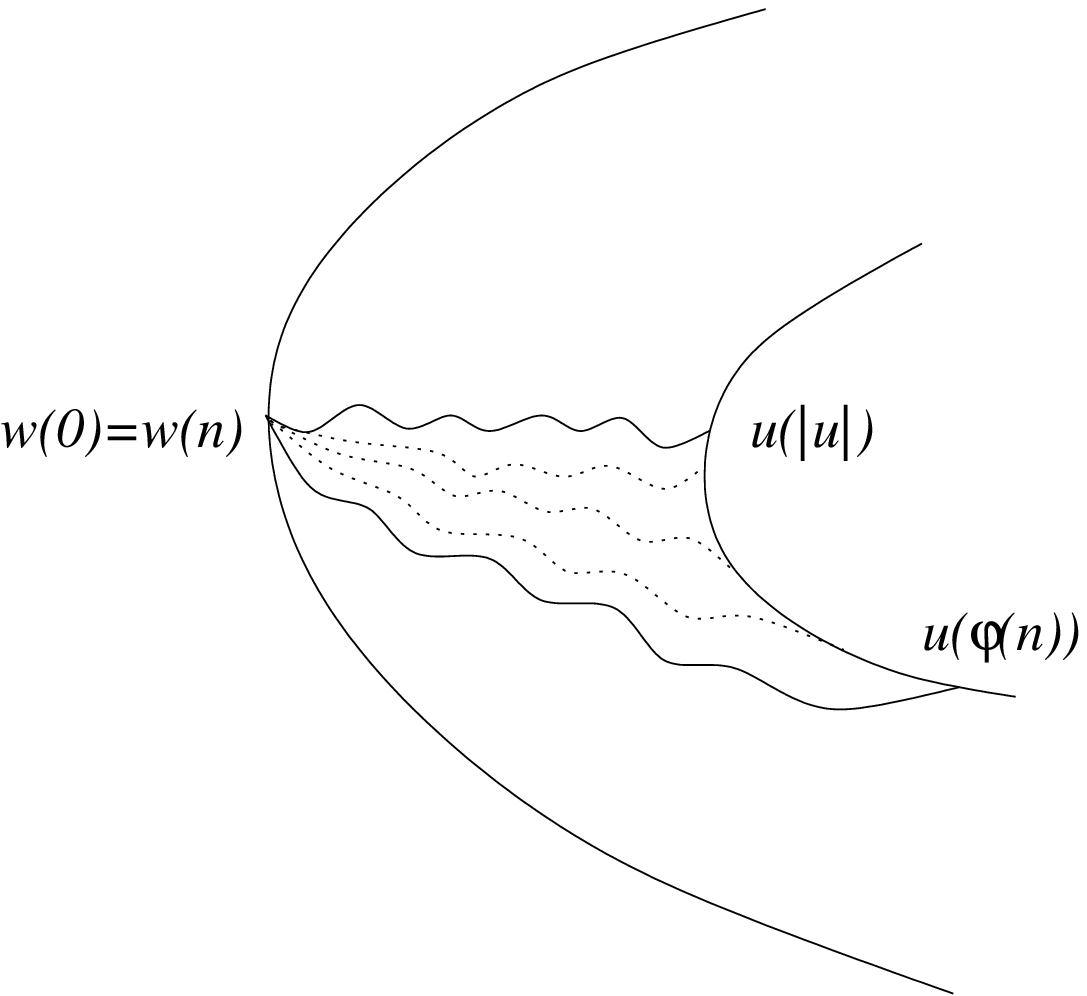}
  \end{center}
  \caption{$\phi(n)<|u|$}
  \label{figasynch1}
\end{figure}

The function $\phi'$ is proper and monotone increasing. In addition,
$w$ and $u$ asynchronously $(k+\epsilon)$-fellow travel with respect
to $\phi'$ since for $ n-\epsilon \leq t \leq n$ we have
$d(w(t),u(\phi'(t)))\leq d(w(t),w(n))+d(w(n),u(\phi'(t)))\leq \epsilon
+k$.

So without loss of generality (by possibly choosing a different
$\phi$) we may assume $n\in E$.  We divide the argument into three
cases.

\noindent
\textbf{Case 1}: If $|t-\phi(t)|\leq 2k$ for all $t\in \mathbb N,
t\leq n$ then $d(w(t),u(t))\leq d(w(t),u(\phi(t)))+d(u(\phi(t)),u(t))
\leq k+2k=3k$ so by Lemma \ref{lem:discrete} $w,u$ synchronously
$(3k+1)$-fellow travel.

\noindent
\textbf{Case 2}: If $t-\phi(t)>2k$ for some $t\in \mathbb N, t\leq n$
then let $j=\min\{t\in \mathbb N: t-\phi(t)>2k\}$.  In particular
$j-0>\phi(j)-\phi(0)+2k$.  Let $l=\max\{l\in \mathbb N: l<j,
j-l>\phi(j)-\phi(l)+2k\}$.  Then $j-(l+1) \not > \phi(j)-\phi(l+1)+2k$
since $l$ was chosen to be maximal, so $j-(l+1) \leq
\phi(j)-\phi(l+1)+2k \leq \phi(j)-\phi(l)+2k $ since $\phi$ is
monotone increasing. Thus we have $j-l=\phi(j)-\phi(l)+2k+\delta$ for
some $\delta \in (0,1]$.

For all $t\in \mathbb N$ with $t\in (l,j)$, we have $t-\phi(t) \leq
2k$ since $t<j$. If $\phi(t)-t \geq 0$ then $j-t>\phi(j)-\phi(t)+2k$
and $l$ is not maximal.  Thus $d(w(t),u(t) \leq
d(w(t),u(\phi(t)))+d(u(\phi(t)),u(t)) \leq k+2k=3k$.

 Let
$w_1=[w(0),w(l)],w_2=[w(l),w(j)],w_3=[w(j),w(n)],u_1=[u(0),u(\phi(l))],
u_2=[u(\phi(l)),u(\phi(j))], u_3=[u(\phi(j)),u(|u|)]$, $p_1$ a path of
length $k_1$ from $w(l)$ to $u(\phi(l))$, and $p_2$ a path of length
$k_2$ from $u(\phi(j))$ to $w(j)$, as in Figure \ref{figasynch2}.
\begin{figure}[ht!]
  \begin{center}
     \includegraphics[width=7cm]{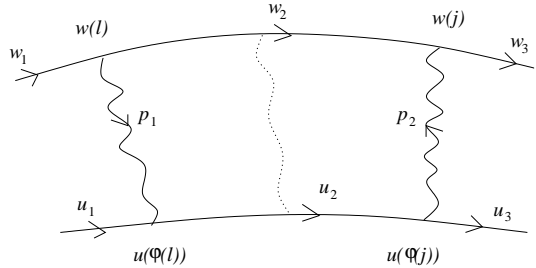}
  \end{center}
  \caption{Case 2: $j-l>\phi(j)-\phi(l)+2k$ }
  \label{figasynch2}
\end{figure}

Define $v=w_1p_1u_2p_2w_3$ which is seen as the bold path in Figure
 \ref{figasynch2_vbold}.  This loop has length
 $l+k_1+(\phi(j)-\phi(l))+k_2+(n-j) =\phi(j)-j+l-\phi(l)+k_1+k_2+n\leq
 \phi(j)-j+l-\phi(l)+2k+n < n$.
\begin{figure}[ht!]
  \begin{center}
     \includegraphics[width=7cm]{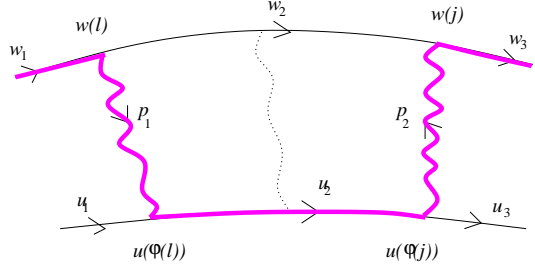}
  \end{center}
  \caption{Case 2: the path $v$ is shown in bold.}
  \label{figasynch2_vbold}
\end{figure}
We will now show that $w$ and $v$ synchronously fellow travel.  For
 $0\leq t\leq l$ the paths $w,v$ $0$-fellow travel.  For $l<t \leq
 l+k_1$ we can find a path of length at most $2k_1$ from $w(t)$ back
 along $w$ to $w(l)$ then down $p_1$ to $v(t)$. Thus $d(w(t),v(t))\leq
 2k_1 \leq 2k$.

The vertex $u(\phi(l))=v(l+k_1)$, so
$d(u(l+k_1),u(\phi(l)))=|(l+k_1)-\phi(l)|\leq 2k+k_1\leq 3k$.  For
$t\in \mathbb N$ and $l+k_1<t\leq l+k_1+\phi(j)-\phi(l)$ we have $t<j$
so $d(u(\phi(t)),u(t))\leq 2k$ and so $d(w(t),v(t))\leq
d(w(t),u(\phi(t)))+d(u(\phi(t)),u(t))+ d(u(t),v(t))\leq k+2k+3k=6k$.

Now $j-\phi(j)=l-\phi(l)+2k+\delta$ so $l-\phi(l)+\phi(j)=j-2k-\delta$
and so $l+k_1-\phi(l)+\phi(j)=j+k_1-2k-\delta$. Let
$r=2k+\delta-k_1\leq 2k+1$.

The vertex $u(\phi(j))=v(j-r)$ so for $t\in \mathbb N$ and  $j-r<t\leq
j-r+k_2$ there is a path from $w(t)$ along $w$ to $w(j)$ then down
$p_2^{-1}$ to $v(t)$ of length at most $r+k_2\leq 3k+1$, so
$d(w(t),v(t))\leq 3k+1$.

Now $d(w(j-r+k_2),w(j))=r-k_2\leq 2k+1$ and for $t \in \mathbb N$ and
$j-r<t\leq n$ $w$ and $v$ travel at constant speed along $w$, at
distance $r-k_2$ apart, so $d(w(t),v(t))\leq 2k+1$.
Thus in total, and by Lemma \ref{lem:discrete}, $w$ and $v$
synchronously $(6k+1)$-fellow travel.

\noindent
\textbf{Case 3}: If $t-\phi(t)\leq 2k$ for all $t\in \mathbb N$ and
$t\leq n$ but $\phi(t)-t> 2k$ for some $t\in \mathbb N$, $t\leq n$,
then let $j=\max\{t\in \mathbb N: t\leq n,\phi(t)-t> 2k\}$.

Then $\phi(j+1)-(j+1)\leq 2k$ so $\phi(j)-j\leq \phi(j+1)-j\leq 2k+1$
so $\phi(j)-j=2k+\delta$ for some $0<\delta\leq 1$.

Let
$w_1=[w(0),w(j)],w_2=[w(j),w(n)],u_1=[u(0),u(\phi(j))],u_2=[u(\phi(j)),u(|u|)]$,
$p_1$ a path of length $k_1$ from $w(j)$ to $u(\phi(j))$, and $p_2$ a
path of length $k_2$ from $u(\phi(|u|))$ to $w(n)$, as in Figure
\ref{figasynch3}.
\begin{figure}[ht!]
  \begin{center}
    \includegraphics[width=5cm]{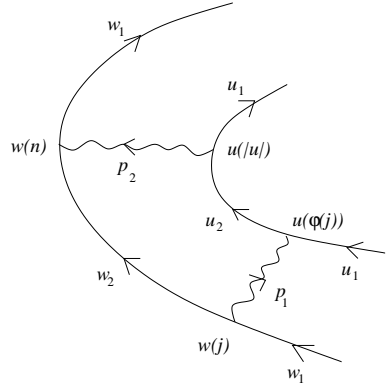}
  \end{center}
  \caption{Case 3: $\phi(j)-j> 2k$}
  \label{figasynch3}
\end{figure}

Define $v=w_1p_1u_2p_2$, shown in bold in Figure
\ref{figasynch3_vbold}.
\begin{figure}[ht!]
  \begin{center}
    \includegraphics[width=5cm]{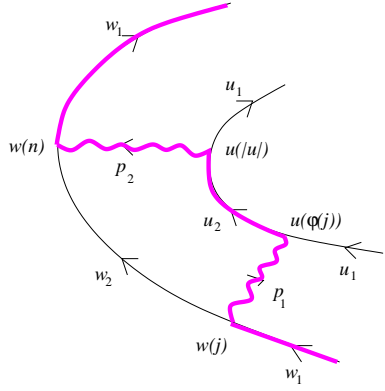}
  \end{center}
  \caption{Case 3: the path $v$ is shown in bold.}
  \label{figasynch3_vbold}
\end{figure}
 This loop has length $j+k_1+(|u|-\phi(j))+k_2 \leq |u| - (\phi(j)-j)
+2k < |u|\leq n$.  We will show that $w$ and $v$ synchronously fellow
travel.

For $0\leq t\leq j$ the paths $w$ and $v$ synchronously $0$-fellow
travel.  For $j< t \leq j+k_1$ there is a path from $w(t)$ back along
$w_2$ to $w(j)$ then down $p_2$ to $v(t)$ of length at most $2k_1$, so
$d(w(x),v(x))\leq 2k_1 \leq 2k$.

The vertex $u(\phi(j))=v(j+k_1)$, so for $t\in \mathbb N$ and
$j+k_1<t\leq j+k_1+|u|-\phi(j)$ we have
$d(u(t),v(t)=\phi-j-k_1=2k+\delta-k_1\leq 2k+1$ so $d(w(t),v(t))\leq
d(w(t),u(\phi(t))) + d(u(\phi(t)),u(t))+ d(u(t),v(t))\leq
k+2k+2k+1=5k+1$.

Recall that $n\in E$ by the argument at the start of this proof, so
$\phi(n)=|u|$.  Now since $n>|u|$ then $n-\phi(n)>0$ so it must be
that $j<n$, and so $n-\phi(n)\leq 2k$.

We have $j+k_1+|u|-\phi(j)=|u|+k_1+j-\phi(j)=\phi(n)+k_1-(2k+\delta)$.
Now for $t\in \mathbb N, j+k_1+|u|-\phi(j) <t \leq
j+k_1+|u|-\phi(j)+k_2$ there is a path from $w(t)$ along $w$ to $w(n)$
of length at most $n-(j+k_1+|u|-\phi(j))$ and from $w(n)$ there is a
path down $p_2^{-1}$ to $v(t)$ of length at most $k_2$. Thus
$d(w(t),v(t))\leq
n-(j+k_1+|u|-\phi(j))+k_2=n-\phi(n)+\phi(j)-j-k_1+k_2=n-\phi(n)+2k+\delta-k_1+k_2\leq
2k+2k+1+k=5k+1$.

Thus in total and by Lemma \ref{lem:discrete} $w$ and $v$
synchronously $(5k+2)$-fellow travel.

Finally for the basepoint case, we merely repeat the argument with
$d(w(0),u(0))=0$.

\end{proof}

\section{Quadratic \iif\ and \ac ity}
In this section we establish connections between the two loop
properties, quadratic \iif s and \ac ity.

\begin{thm}
\label{loopquad}
If $(G,X)$ has the \lsp\ then $G$ is finitely presented, and has
a quadratic \iif.
\end{thm}

\begin{proof}
Let $w=_G1$. Define $w=w_0$.  While $w_i$ is not the empty word, there
 is a shorter loop $w_{i+1}$ that $k$-fellow travels $w_i$.  After at
 most $|w|$ iterations we get the trivial word. The space between
 $w_i$ and $w_{i+1}$ can be filled by $|w_i|$ relations of length at
 most $2k+2$, so it follows that $G$ is finitely presented as $\langle
 X | \{u\in X^*:|u|\leq 2k+2\}\rangle$.  Moreover, the number of such
 relations needed to fill $w$ is at most $\Sigma_{i=1}^{|w|} i \leq
 |w|^2$.
\end{proof}

\begin{thm}
\label{bploopac}
If $(G,X)$ has the basepoint \lsp\ then $(G,X)$ \ac.
\end{thm}

\begin{proof}
Let $w$ and $u$ be two geodesics of length $N$ such that $d(\overline
w, \overline u)\leq 2$, realized by a path $\gamma$.  Let $k$ be the
basepoint loop shortening constant, and without loss of generality
assume it is an even integer.  The word $w\gamma u^{-1}$ is a loop
based at the identity vertex, of length at most $2N+2$.  Applying
basepoint loop shortening we get a loop $y$ based at the identity of
length at most $2N+1$. Applying the property once more we get a loop
$v$ based at the identity of length at most $2N$, so $v \subseteq
B(N)$.

The path that retraces $w$ back to $w(N-\frac{k}{2})$ , then travels
across to $y(N-\frac{k}{2})$ then to $v(N-\frac{k}{2})$, then travels
along $v$ to $v(N+\frac{k}{2}+2)$, over to $y(N+\frac{k}{2}+2)$ then
to $u(N+\frac{k}{2}+2)$, then along $u$ to its end lies in $B(N)$ and
has length at most $6k+2$. See Figure \ref{fig:acproof}.
\begin{figure}[ht!]
  \begin{center}
     \includegraphics[width=10cm]{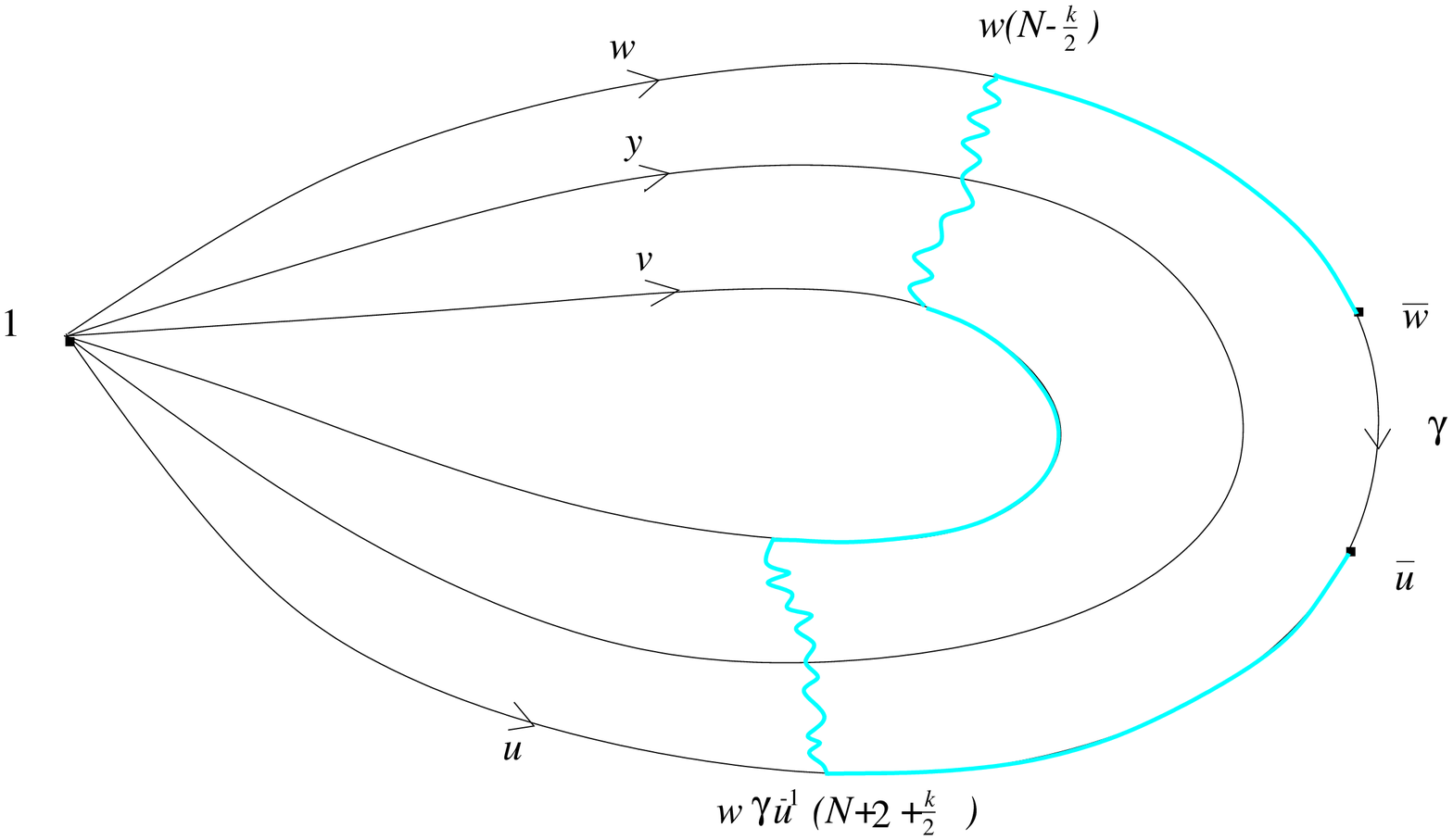}
  \end{center}
  \caption{The basepoint \lsp\ implies \ac}
  \label{fig:acproof}
\end{figure}
\end{proof}

The theorem provides an easy route to proving \ac\ for some examples,
and is potentially an extremely useful tool.
In the next section we show this by ``reproving'' a theorem of the
author in \cite{\Enonhopf}.

\begin{figure}[ht!]
  \begin{center}
      \includegraphics[width=11cm]{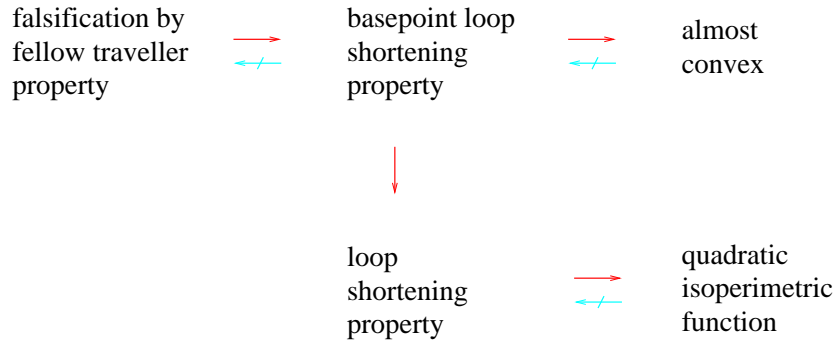}
  \end{center}
  \caption{Implication diagram}
  \label{fig:implic}
\end{figure}

We summarize the results so far in Figure \ref{fig:implic}.
The non-reversible implications (in grey
)
will be proved by counterexamples below.

\section{Multiple \hnn s}
In \cite{\Enonhopf} the author proves that a certain class of multiple
\hnn\ group presentations are \ac. We show that the same hypothesis
implies the basepoint \lsp.  This gives a rapid proof of \ac ity for
some examples of interest.

\begin{defn}[Multiple \hnn] \label{mult}
Let $(A,Z)$ be a group with finite \gset\ $Z$ and relations $R$,  let
 $U_1,\ldots, U_n, V_1, \ldots, V_n$ be subgroups of $A$ and let
$\phi_i:U_i \ra V_i$ be an isomorphism for each $i$.
The group $(G,X)$ with presentation
$$\langle Z,s_1,\ldots, s_n | R, s_i^{-1}u_is_i=\phi_i(u_i) \; 
\forall u_i\in U_i, \; \forall i \rangle$$
is a {\em multiple \hnn} of $(A,Z)$. The generators $s_i$ are called
{\em stable letters}, 
and the pairs of $U_i,V_i$ are called {\em associated subgroups}.
\end{defn}

\noindent
If each $U_i$ is finitely generated by $\{ u_{i_j} \}$ and $\phi_i (u_{i_j})=v_{i_j}$
then $V_i$ is finitely generated by $\{v_{i_j}\}$.
Thus $(G,X)$ has the finite presentation
$$\langle Z,s_1,\ldots, s_n | R, s_i^{-1}u_{i_j}s_i=v_{i_j} \; 
\forall i, \; \forall j \rangle.$$

\begin{thm}[Britton's Lemma]
Let $(G,X)$ be a multiple \hnn\ with the presentation in Definition
\ref{mult} above. If $w\in X^*$
is freely reduced and $w=_G1$ then $w$ contains a sub-word of the form
$s_i^{-1}u_{i_j}s_i$ or $s_iv_{i_j}s_i^{-1}$ for some non-trivial
$u_{i_j}\in U_i$ or $ v_{i_j}\in V_i$.
\end{thm}

\noindent
A sub-word $s_i^{-1}u_{i_j}s_i$ or $s_iv_{i_j}s_i^{-1}$
is called a {\em pinch}, and if a word admits no pinches it is called
{\em stable letter reduced}.

\begin{defn}[Strip equidistant]
Let $(G,X)$ be a multiple \hnn\ with the presentation in Definition
\ref{mult} above.
If $|u_i|=|\phi(u_i)|$ for all $i$ then we say $(G,X)$ has a 
{\em strip equidistant} presentation.
\end{defn}

\noindent
Note that if $(G,X)$ has a strip equidistant presentation then if a word $w\in X^*$
admits a pinch then it can be shortened by $2$, so geodesics are stable letter reduced.


\begin{defn}[Totally geodesic]
Let $(G,X)$ be any group with \gset\ $X$. 
A subgroup $A$ of $G$ with \gset\ $Y\subseteq X^*$ 
 is {\em totally geodesic} in $(G,X)$ if every geodesic word $w\in X^*$ 
evaluating to an element of $A$  is an element of $Y^*$.
\end{defn}


\begin{thm}
\label{multhnnbl}
Let $(G,X)$ be a multiple \hnn\ of $(A,Z)$ as in Definition \ref{mult} with
a strip equidistant presentation, such that
 associated subgroups are totally geodesic and
 $(A,Z)$ enjoys the \fftp.
Then $(G,X)$ enjoys the basepoint \lsp.
\end{thm}

\begin{proof}
Let $k$ be the
\fftp\ constant for $(A,X)$.
Let $w$ be a loop based at $w(0)$ in $(G,X)$.
If $w$ has no stable letters, then there is a shorter loop in
$(A,X)$ that $k$-fellow travels $w$ by the \fftp\ in $(A,X)$.

If $w$ has stable letters, then by Britton's Lemma it admits a pinch.
Let $sw_2s^{-1}$ be an inner-most pinch, that is, $w_2 \in Z^*$, and
$w=w_1sw_2s^{-1}w_3$.
If $w_2$ is not a geodesic then apply the \fftp\ in $(A,X)$ to  get a
shorter sub-word $u$ so that $w_1sus^{-1}w_3$ $k$-fellow travels $w$.
If $w_2$ is geodesic then by total geodecity of associated subgroups, 
$w_2$ is a word in $\{u_i\}^*$ [respectively  $\{v_i\}^*$].
Then $sw_2s^{-1}$ is 2-fellow traveled by $\phi (w_2)$ 
[respectively  $\phi^{-1}(w_2)$].
\end{proof}

One might think that the preceding proof
can be strengthened to show that 
$(G,X)$ in fact has the \fftp. The first example of the next section 
shows that this is not possible. Moreover the last example of the next
section shows that the totally geodesic hypothesis
cannot be relaxed.

\section{Examples}
In this section we consider four group presentations which display a
diverse range of  properties. In particular we will fill in
 Table \ref{table} below of examples and their properties.
\begin{table}[ht!]
\caption{Examples and their properties}
\label{table}
\[
\begin{tabular}{|l||c|c|c|c||c|c|c|}
\hline
 & falsification & basepoint & loop      & almost & automatic & \cat(0) & quadratic\\
 & by fellow     & loop      & shortening& convex  &        &   & Dehn \\
 & traveler      & shortening& property  &   &       &   & function\\
 & property      & property  &           &    &      &   &  \\
\hline\hline
$(G_B,X)$ & ?  & yes  & yes  & yes & yes  & yes  & yes \\
\hline
$(G_W,X)$ & no  & yes  & yes  & yes & ?  & yes  & yes \\
\hline
$(G_G,X)$ & no  & no  & no   & yes &  no  & no   & yes  \\
\hline
$(G_S,X)$ & no  & no  & ?   & no  &  no   & no   &  yes  \\
\hline
\end{tabular}
\]
\end{table}

\begin{eg}\label{wiseeg}
$(G_W,X)=\langle a,b,c,d,s,t|c=ab, c=ba, d=c^2, s^{-1}as=d,
t^{-1}bt=d \rangle$.
\end{eg}

\noindent
This group  was  considered by Wise \cite{\Wise}, who proved it is
\cat(0) and non-Hopfian. The group is a double \hnn\ of
$(\Z,\{a,b,c,d\})$, with associated subgroups $\langle a \rangle,
 \langle b \rangle, \langle d \rangle$  totally geodesic in $(\Z^2, \{a,b,c,d\})$.
Neumann and Shapiro prove that any finite \gset\ for an abelian group
has the \fftp\ \cite{\NSgeomfinite}.
 It follows from Theorem \ref{multhnnbl} that $(G_W, X)$
enjoys the basepoint \lsp\ and is consequently \ac.
 The author has shown that this example
does not enjoy the \fftp\ (see \cite{\Ethesis}), and so the basepoint
\lsp\ does not imply the \fftp.
It is not known whether this group is automatic.

\begin{eg}
\label{bridsoncat0eg}
$(G_B,X)=\langle a,b,\gamma, s,t| c=aba^{-1}b^{-1},
\gamma a \gamma ^{-1}=a^{-1}, 
\gamma b \gamma ^{-1}=b^{-1},$

\noindent
$sas^{-1} =c, tbt^{-1}=c \rangle$.
\end{eg}

\noindent
Bridson shows that this group cannot act on any 2-dimensional \cat(0)
space,
but is the fundamental group of a 3-dimensional non-positively curved
cube complex \cite{\Bridcat}. It follows from Niblo and Reeves
\cite{\NibloR} that $G_B$ is biautomatic.
The pair is a triple \hnn\ of $F_2$, the free group on two letters.
Since $F_2$ is word-hyperbolic, it enjoys the \fftp\ with respect to
any \gset.
It is easy to see that the presentation is strip equidistant, and 
the associated subgroups are totally geodesic, and so this pair
has the basepoint \lsp\ by Theorem  \ref{multhnnbl} and consequently is \ac.
It is  not known whether it enjoys the \fftp.

One might ask whether  every group with quadratic \iif\ is
\ac. The following example shows this is not the case, and in addition
shows that the basepoint loop property is not equivalent to the
enjoyment of a quadratic \iif.

\begin{eg}
\label{stall}
Let $\phi:(F_2)^3\rightarrow \Z$ be a homomorphism which sends each
word in $(F_2)^3=\langle a,b,c,d,e,f \rangle$ to its exponent sum.
Define the group $G_S=Ker(\phi)$.
\end{eg}

\noindent
Stallings showed that $G_S$ is finitely presented but is not of
type FP$_3$ \cite{\Stallings}, and so not of type F$_3$.
Recently Bridson has shown this group has a quadratic \iif\
\cite{\Bridson}.
Since $G_S$ is not of type F$_3$, it is not automatic, not \cat(0),
and does not enjoy the \fftp\ for any \gset\ (See \cite{\GerstenGGT,\BridsonHaefliger,\Efftpft} respectively).
Finding a \gset\ for which  $G_S$ is \ac\ or has the [basepoint]
\lsp\ would prove that the respective property does not imply F$_3$.

\begin{figure}[ht!]
  \begin{center}
      \includegraphics[width=3.5cm, height=4.5cm]{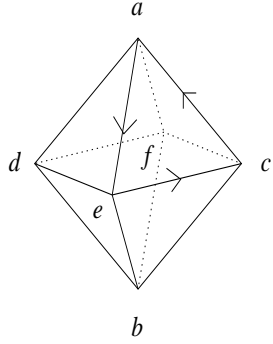}
  \end{center}
  \caption{Flag complex encoding $(G_S,X)$}
  \label{fig:octahedron}
\end{figure}

There is a standard way of associating a right-angled Artin group to a 
finite flag complex (See  Bestvina and Brady \cite{\BestBrady} for
details and references).
Dicks and Leary give the following description of a presentation for
$G_S$, based on work of Bestvina and Brady \cite{\BestBrady}.
Consider an octahedron with opposite vertices labeled by generators of
each free factor of $(F_2)^3$. See Figure \ref{fig:octahedron}.
 Each directed edge of the octahedron defines a generator of $G_S$, where
$ac^{-1}$ is the edge from $c$ to $a$. For convenience we 
denote the inverse of a generator of $(F_2)^3$ in upper case, so the
edge $ac^{-1}$ is written $aC$, and so on.
Each 2-cell of the octahedron defines two relations, so the
2-cell with vertices $a,c,e$ defines two relations $aCcEeA$ and $aCeAcE$.
Dicks and Leary prove in \cite{\DicksLeary} that these 12 generators
and 16 relations are a presentation for $G_S$. We will denote this
\gset\ by $X$.

Define a homomorphism
$\rho:G_S \rightarrow (F_2)^3$ by $\rho(yZ)=yZ$.
It is clear that if two words $u,v$ evaluate to the same element of
$G_S$ then $\rho(u)=\rho(v)$ in $(F_2)^3$.
\begin{lem}
\label{BB}
Let $w\in X^*$ and $z \in \{c,d,e,f\}$.
If $\rho(w)= f^nd^nE^nC^{n-1}Z$ in $(F_2)^3$
then $|w|\geq 3n$.
\end{lem}
\begin{proof}
The first $E$ or $C$ in $w$ that does not freely cancel with an
earlier letter must occur after $w(n)$, 
since $f^n$ or $d^n$ must be read before 
$E$ or $C$ respectively.
After $w(n)$ we must read $E^nC^{n-1}Z$ so we need at 
least $2n$ letters of $X$. So $|w|\geq 3n$.
\end{proof}
It follows that $\alpha =(fA)^n(dE)^n(aC)^{n-1} , \beta =
 (fB)^n(dE)^n(bC)^{n-1}$
 are geodesic since they are sub-words of $\alpha aC, \beta bC$ 
which are geodesic by the lemma.
\begin{lem}
\label{BB2}
Let $w\in X^*$, $z \in \{c,d,e,f\}$, and $v\in X^*$
 such that $\rho(v)= 1$ in $\langle a,b \; | \; - \rangle $.
If $v$ is of  length less than $n-1$ and $\rho(w)= f^nd^nE^nC^{n-1}vZ$ in $(F_2)^3$
then $|w|\geq 3n$.
\end{lem}
\begin{proof}
$v$ is  shorter than $n-1$, so cannot freely cancel all the $E^n$ and $C^{n-1}$.
So again the first $E$ or $C$ in $w$ must occur after $w(n)$.
If $v$ freely cancels some of the $E^n, C^{n-1}$ (and $Z$) $v$ must
 contain $e^i,c^j$ (and $z$)
so must contain $i+j (+1)$ upper case letters, since $v \in X$.
In sum total $E^nC^{n-1}vZ$ has at least $2n$ upper case letters 
(plus possibly more upper and lower pairs).
So again there must be at least $2n$ letters
in $X$ after $w(n)$, so $|w|\geq 3n$.
\end{proof}

\begin{thm}
$(G_S,X)$ is not \ac.
\end{thm}
\begin{proof}
 Let  $\alpha =(fA)^n(dE)^n(aC)^{n-1} , \beta = (fB)^n(dE)^n(bC)^{n-1}, \gam =aCcB$.
It is easily checked using Lemma \ref{BB}  that $\alpha, \beta$ are geodesics, each
of length $3n-1$ ending distance 2 apart in $\Gamma(G_S,X)$, realized by $\gamma$.
Assume by way of contradiction that $(G_S,X)$ is \ac,
 so there is a path $p$ from $\overline{\alpha} $ to
 $\overline{\beta} $ inside $B(3n-1)$ of bounded length.
 See Figure \ref{fig:4_14}.
\begin{figure}[ht!]
  \begin{center}
      \includegraphics[width=8cm]{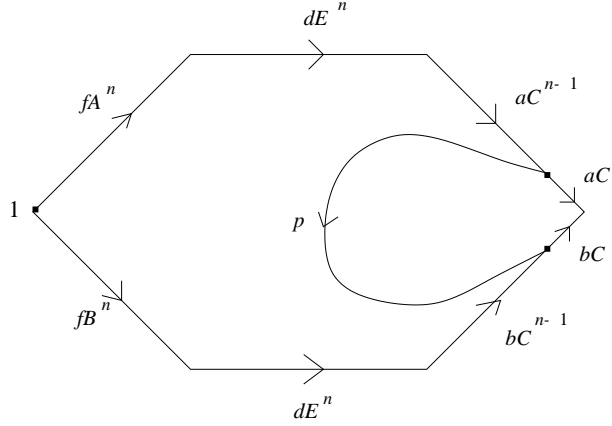}
  \end{center}
  \caption{Paths in $\Gamma(G_S,X)$}
  \label{fig:4_14}
\end{figure}
We have $pbCcA=_{G_S} 1$ so
$pbCcA=_{(F_2)^3} pbA=_{(F_2)^3} 1$.
This means $p$ must contain a $B$ to cancel with the $b$ in this word, so $p=uzBv$
with $z \in \{c,d,e,f\}$ and $\rho(v)= 1$ in $\langle a,b \; | \; - \rangle $, having bounded length since 
$p$ is of bounded length.
We choose $n$ to be greater than this bound.
Let $g$ be a geodesic to $\overline{\alpha u}$, as in Figure \ref{fig:4_15}. 
\begin{figure}[ht!]
  \begin{center}
      \includegraphics[width=8cm]{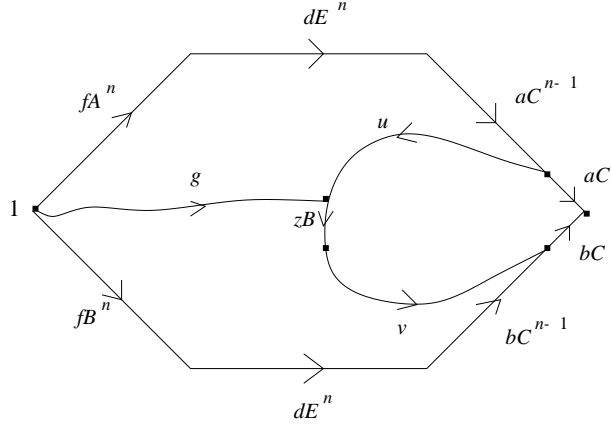}
  \end{center}
  \caption{A geodesic $g$ to $\overline{\alpha u}$}
  \label{fig:4_15}
\end{figure}
Now
$\rho(g) = (fB)^n(dE)^n(bC)^{n-1}v^{-1}bZ  
 = f^nd^nE^nC^{n-1}(v^{-1})Z$
in $(F_2)^3$.
By Lemma \ref{BB2} we have $|g|\geq 3n$, which contradicts the fact that 
$p\subseteq B(3n-1)$.
\end{proof}

\begin{figure}[ht!]
  \begin{center}
     \includegraphics[width=12.5cm]{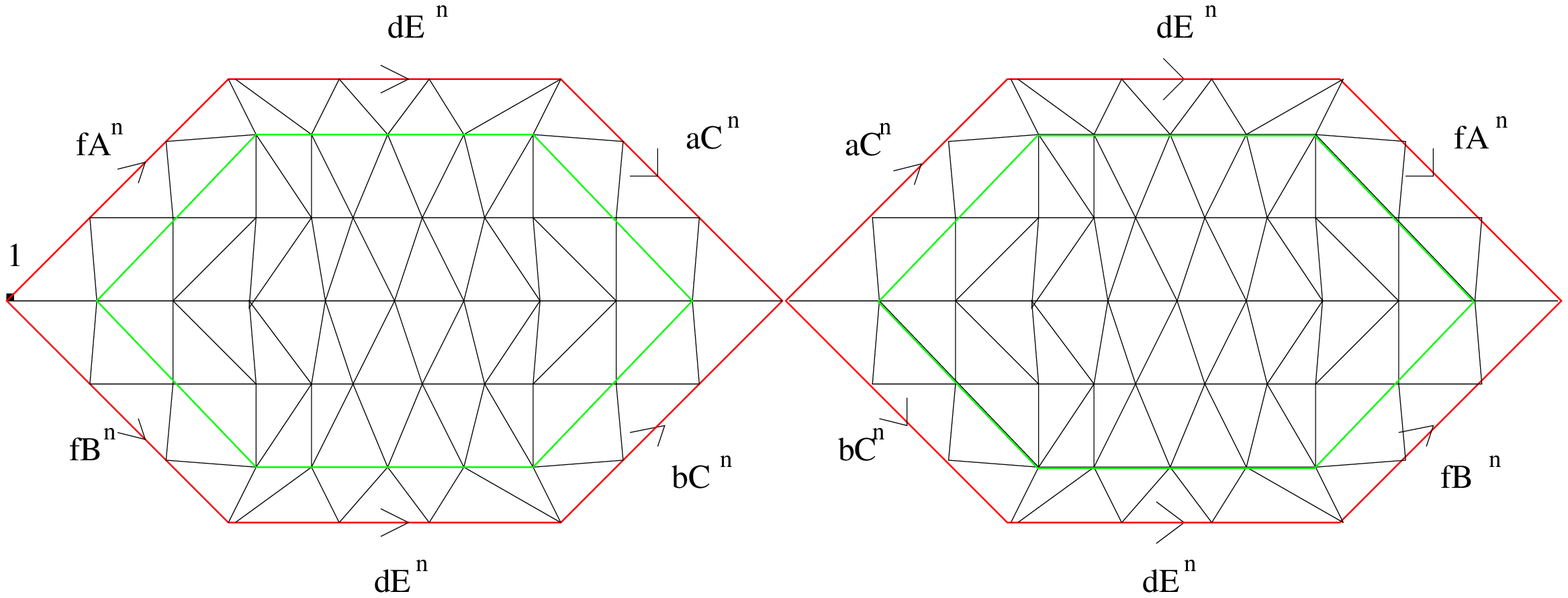}
  \end{center}
  \caption{A potential counterexample to the \lsp\ for $(G_S,X)$}
  \label{badloopStallBB}
\end{figure}

The author has considered alternate finite presentations for $G_S$;
see \cite{\Ethesis}. It may be that this example is \ac\ for another
(possibly weighted) \gset.
The boundary loop shown in Figure \ref{badloopStallBB}
does not appear to be fellow traveled by a shorter loop
for any constant independent of $n$. It is likely (but not proved)
that $(G_S,X)$ does not have the \lsp.

The final example is another multiple \hnn, but does not satisfy the
totally geodesic associated subgroup hypothesis of Theorem \ref{multhnnbl}.
It is \ac\ however, 
and has a quadratic \iif, so one might suspect
that it would have the loop shortening property.

\begin{eg}
\label{gersteg}
$(G_G,X)=\langle a,b,c,d,s,t|c=ab, c=ba, d=ab^{-1},  s^{-1}as=c,
t^{-1}at=d \rangle$.
\end{eg}

\noindent
Gersten proves that this group is not \cat(0) \cite{\Gerstencat}.
Brady and Bridson
showed that the group is has quadratic \iif\ \cite{\BBIsop} 
and is not biautomatic \cite{\BB}. 
The author proves the pair is \ac\ and fails the \fftp\ in \cite{\Epatterns}.        
Recent work of Bridson and Reeves shows that $G_4$ is not automatic.
The group is free-by-cyclic.

In  \cite{\Epatterns}  the \cg\ of $(G_G,X)$ is described as being made
up of copies of the \cg\ for $(\Z^2,\{a,b,c,d\})$, which we call
``planes'',
glued together along bi-infinite lines $a^i,c^i,d^i$ by stable letter
``strips''.
We now give a more technical definition of the idea of a strip.
\begin{defn}
An {\em $s$-strip} is the set of open edges of the form 
$\{(wa^i,wa^is):i\in \Z\}$
 for  some arbitrary word $w$.
We denote this strip by $(w\langle a\rangle,s)$.
The three other  possible strips  are
 $(w\langle a\rangle,t)$, $(w\langle c\rangle,s^{-1})$ $(w\langle d\rangle,t^{-1})$.
\end{defn}

\begin{lem}
A strip divides the \cg\ into two connected half spaces.
\end{lem}
\begin{proof}
Let $(w\langle x\rangle,r)$ be a strip in $\Gamma(G_G,X)$.
Since $(G_4,X)$ is strip equidistant, a geodesic crosses each strip at
most once.
Let $H_-$ be the set of all points in  $\Gamma(G_G,X)$
so that a geodesic from it to $w$ does not cross the strip.
Let $H_+$ be the set of all points in  $\Gamma(G_G,X)$
so that a geodesic from it to $wr$ does not cross the strip.
It is easily seen that the \cg\ is the (disjoint) union of
$H_-, H_+$ and the strip $(w\langle x\rangle,r)$, 
 the two components
are each path connected, and  $H_-\cap H_+=\emptyset$.
\end{proof}
As a consequence we can say that two points lie on the {\em same side}
of a strip if they lie in the same half space.
\begin{thm}
\label{BBnotlsp}
 $(G_G,X)$ does not enjoy the \lsp.
\end{thm}

\begin{proof}
Assume by way of contradiction that $(G_G,X)$ has the \lsp\ with
constant $k$.
Let  $w=d^nst^{-1}c^nd^nst^{-1}c^{-n}d^{-n}st^{-1}c^{-n}d^{-n}st^{-1}c^n$ for
$n>k$.
See Figure \ref{fig:G_4loop}.
\begin{figure}[ht!]
  \begin{center}
      \includegraphics[width=11cm]{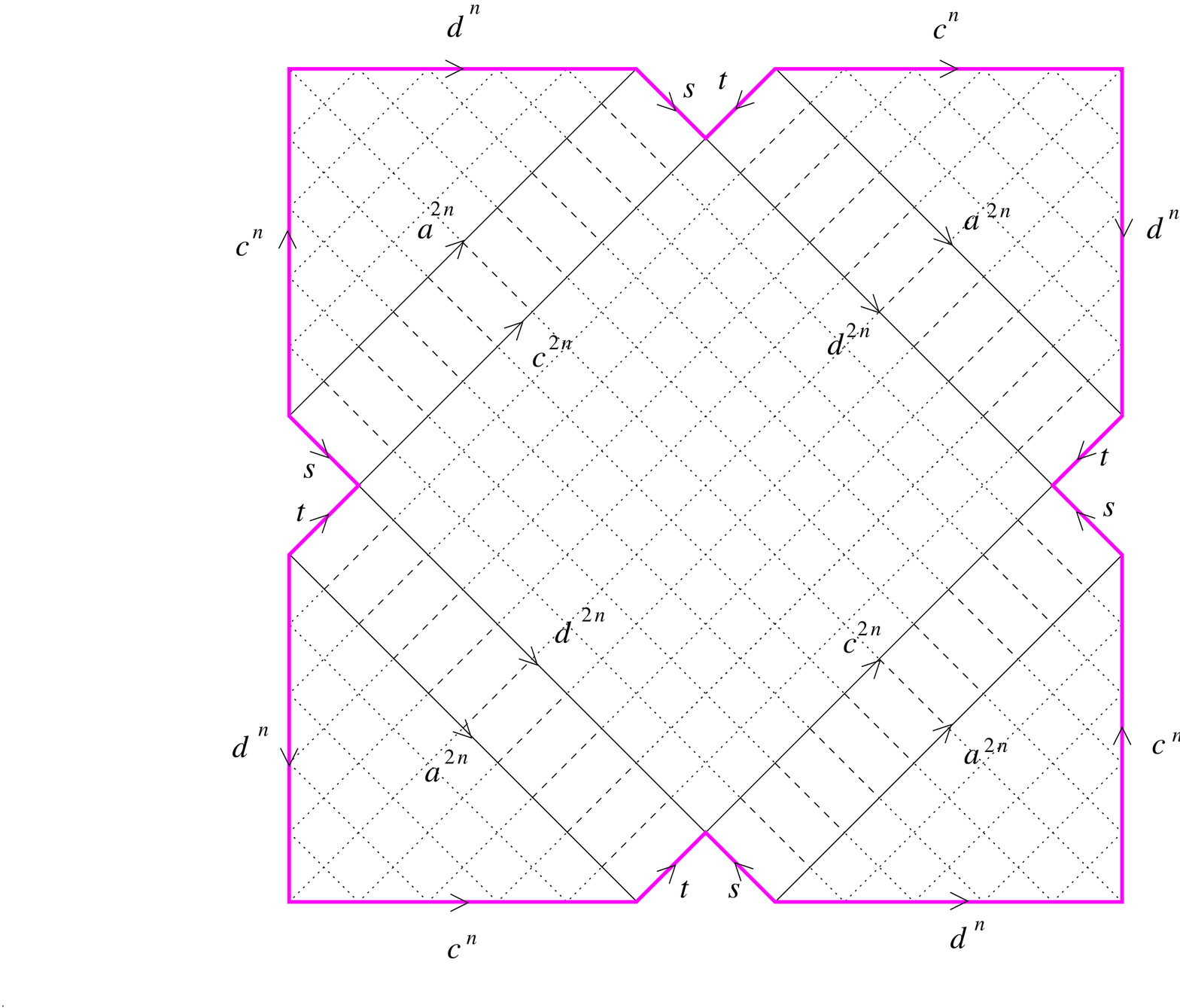}
  \end{center}
  \caption{The loop $w$ in $(G_G,X)$}
  \label{fig:G_4loop}
\end{figure}
It is easy to check algebraically that $w$ is a loop.

Now by assumption there is a loop $u$ of length $|u|<8n+8$ that
synchronously $k$-fellow travels $w$.
The  point $u(0)$ lies on the same side of the strip $S_1$
as $w(0)$. To see this, let $g$ be a geodesic from 
$w(0)$ to $u(0)$. If they lie in different half-spaces, let $(p,ps)$
be the first edge that $g$ crosses on the strip $S_1$.
Then $|g|=d(w(0),p)+1+d(ps,u(0))\geq n+1$.
This is a contradiction since $k<n$.

Repeating the argument, we have that $u(2n+2),u(4n+4),u(6n+6)$ lie on the
same side of the strips $S_2,S_3,S_4$ as $w(2n+2),w(4n+4),w(6n+6)$
respectively.

Now the path $u$ must go between these four points by passing through
the base plane, shown in Figure \ref{fig:baseplaneG_4}.
\begin{figure}[ht!]
  \begin{center}
       \includegraphics[width=11cm]{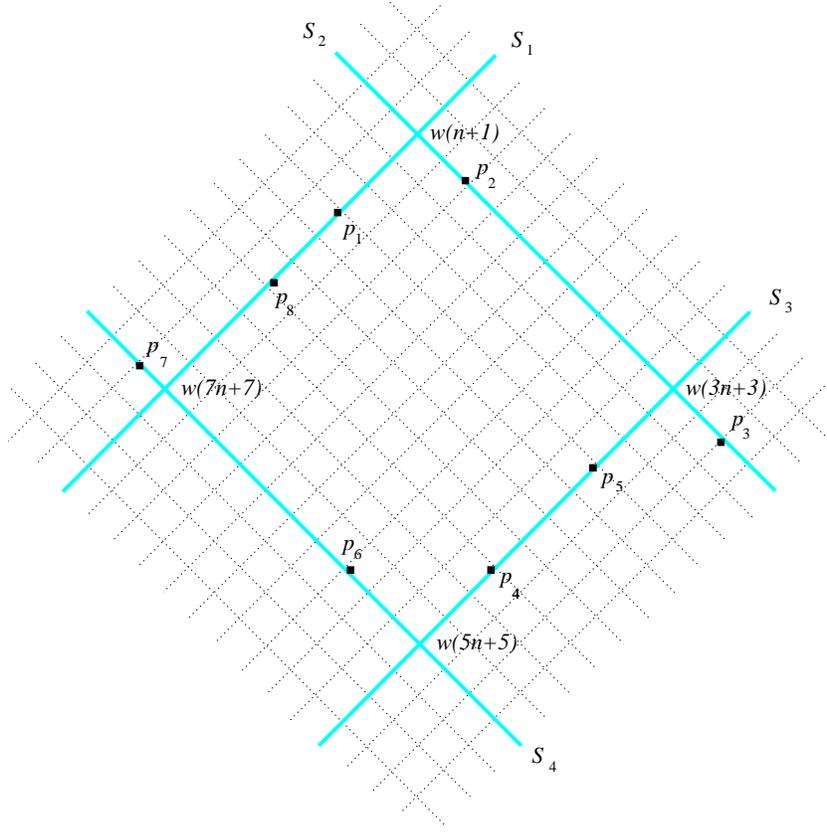}
  \end{center}
  \caption{The ``base plane''}
  \label{fig:baseplaneG_4}
\end{figure}
Let $p_1$ be the first point on the base plane that $u$ crosses after
$u(0)$, $p_2$ the last point on the base plane before $u(2n+2)$,
$p_3$ be the first point on the base plane that $u$ crosses after
$u(2n+2)$, $p_4$ the last point on the base plane before $u(4n+4)$,
and so on up to $p_8$.

Let $m_1$ be the distance between $w(n+1)$ and $p_1$,
$m_2$ be the distance between $w(n+1)$ and $p_2$,
$m_3$ be the distance between $w(3n+3)$ and $p_3$,
$m_4$ be the distance between $w(3n+3)$ and $p_4$,
and so on up to $m_8$.

Notice that we impose no restrictions on where the path $u$ enters and
exits the plane, just that it does so at least eight times, via the
appropriate strips.

Now $d(p_i,p_{i+1})=m_i+m_{i+1}$ for $i=1,3,5,7$ since
$c^id^j$ is a geodesic in the base plane.

The path $u$ must go from $p_2$ to $p_3$ via $u(2n+2)$, so it must
cross the strip $S_2$.
The distances $m_i$ are the same on either side of the strip.
That is, $d(p_2t^{-1}, w(n))= m_2$ and so on.
Then
 $d(p_2t^{-1},p_3t^{-1})\geq |2n-m_2-m_3|$,
 $d(p_4s^{-1},p_5s^{-1})\geq |2n-m_4-m_5|$, 
 $d(p_6t^{-1},p_7t^{-1})\geq |2n-m_6-m_7|$ 
 and
 $d(p_8s^{-1},p_1s^{-1})\geq |2n-m_8-m_1|$.

\noindent
Thus 
$|u| \geq (m_1+m_2) +(m_3+m_4) +(m_5+m_6)+ (m_7+m_8)$

\noindent
$+8+|2n-m_2-m_3|+|2n-m_4-m_5|+|2n-m_6-m_7|+|2n-m_8-m_1|$

\noindent
$=8+ |2n-(m_2+m_3)| + (m_2+m_3)
+|2n-(m_4+m_5)| + (m_4+m_5)$

\noindent
$+|2n-(m_6+m_7)| + (m_6+m_7)
+|2n-(m_8+m_1)| + (m_8+m_1)$

\noindent
$\geq 8+8n $ 
and this contradicts the fact that $u$ must be shorter that $w$.
\end{proof}
Perhaps $G_G$ has the \lsp\ for another \gset.

\section{Open questions}
The three question marks in Table \ref{table} are open. 
We have seen that the loop properties are closely related to the \gset\ dependent
properties of \ac ity and the \fftp, so it is possible that this
unfortunate family trait is inherited. Can we find an example of a
group that enjoys the 
 [basepoint] \lsp\ with respect to one  \gset\ and not another?
Also, is there an example of a group presentation that has the \lsp\
but not the basepoint \lsp?

\bibliography{refs}
\bibliographystyle{plain}

\end{document}